\documentclass[12pt]{article}
\usepackage{amssymb}
\begin{document}

\def\kbar{{\mathchar'26\mkern-9muk}}  
\def\bra#1{\langle #1 \vert}
\def\ket#1{\vert #1 \rangle}
\def\vev#1{\langle #1 \rangle}
\def\tr{\mbox{Tr}\,}
\def\ad{\mbox{ad}\,}
\def\ker{\mbox{Ker}\,}
\def\im{\mbox{Im}\,}
\def\der{\mbox{Der}\,}
\def\ad{\mbox{ad}\,}
\def\b#1{{\mathbb #1}}
\def\c#1{{\cal #1}}
\def\pt{\partial_t}
\def\px{\partial_1}
\def\bpx{\bar\partial_1}
\def\la{\langle}
\def\ra{\rangle}
\def\nn{\nonumber \\}
\def\etal{{\it et al.}\ }

\newcommand{\sect}[1]{\setcounter{equation}{0}\section{#1}}
\renewcommand{\theequation}{\thesection.\arabic{equation}}
\newcommand{\subsect}[1]{\setcounter{equation}{0}\subsection{#1}}
\renewcommand{\theequation}{\thesection.\arabic{equation}}

\newcommand{\be}{\begin{equation}}
\newcommand{\ee}{\end{equation}}
\newcommand{\ba}{\begin{eqnarray}}
\newcommand{\ea}{\end{eqnarray}}

\title{A Calculus Based on a $q$-deformed Heisenberg Algebra}

\author{B.L. Cerchiai$\strut^{1,2}$ \,
        R. Hinterding$\strut^{1,2}$ \\ J. Madore$\strut^{2,3}$, 
        and \ J. Wess$\strut ^{1,2}$\\\\
        \and
        $\strut^1$Sektion Physik, Ludwig-Maximilian Universit\"at,\\
        Theresienstra\ss e 37, D-80333 M\"unchen
        \and
        $\strut^2$Max-Planck-Institut f\"ur Physik\\
        (Werner-Heisenberg-Institut)\\
        F\"ohringer Ring 6, D-80805 M\"unchen
        \and
        $\strut^3$Laboratoire de Physique Th\'eorique et Hautes Energies\\
        Universit\'e de Paris-Sud, B\^atiment 211, F-91405 Orsay
        }
\date{}

\maketitle

\abstract{We show how one can construct a differential calculus
over an algebra where position variables $x$ and momentum variables
$p$ have be defined.  As the simplest example we consider the
one-dimensional $q$-deformed Heisenberg algebra. This algebra has a
subalgebra generated by $x$ and its inverse which we call the
coordinate algebra. A physical field is considered to be an element of
the completion of this  algebra.  We can construct a derivative which 
leaves invariant  the coordinate algebra and so takes physical fields 
into physical fields.  A generalized Leibniz
rule for this algebra can be found.  Based on this derivative
differential forms and an exterior differential calculus can be
constructed.}

\vfill
\noindent
LMU-TPW 98-12\\
MPI-PhT 98-60\\
\noindent
July 1998\\
\medskip
\eject

\parskip 4pt plus2pt minus2pt

\sect{Introduction}

We would like to show how one can construct a differential calculus
over an algebra where position variables $x$ and momentum variables
$p$ have been defined.  As the simplest example we consider the
one-dimensional $q$-deformed Heisenberg algebra~\cite{Sch92}. This
algebra has a subalgebra generated by $x$ and its inverse which we
call the coordinate algebra. A physical field is defined to be an
element of the completion of the  algebra.  We can construct a
derivative which leaves invariant the elements of the coordinate
algebra and so takes physical fields into physical fields.  A
generalized Leibniz rule for this algebra can be found.  Using this
derivative differential forms and an exterior differential calculus
can be constructed. This is done in Section~4.

On the coordinate algebra it is possible to define an integral purely
algebraically, as the inverse image of the derivative; to a given
function we associate another function as the integral .  We study the
above definition of the integral in an explicit representation of the
algebra, a large class of which are known~\cite{Heb94}. This leads us
in a natural way, in Section~5, to a trace formula for the integral
which produces the well known Jackson integral~\cite{Gas90}. We find
that a form of Stokes' theorem can be proven.  There is an interesting
fact about the summation in the formula for the Jackson integral: even
and odd lattice points must be treated separately. The integral
separates therefore into a sum over even lattice points and a sum over
odd lattice points. Such a phenomenon is not unusual for the
$q$-deformed Heisenberg algebra~\cite{Heb94}. The separation can be
traced to the fact that the derivative is actually a second-order
differential operator~\cite{Sch98}.  With the integral we define in
Section~6 an inner product and therefore a Hilbert space $L_q^2$ and
so the notion of an hermitian and a self-adjoint operator becomes
meaningful. Differentiation becomes a linear operator on $L_q^2$. It
is not hermitian but its square has a self-adjoint extension.  In
Section~7 we introduce a basis for $L_q^2$ in terms of $q$-cosine and
$q$-sine functions~\cite{Koo92}.  This opens the way for quantum
mechanics in the form of wave mechanics. In Section~8 a Schr\"odinger
equation is defined. This allows a probability interpretation with a
probability density and a probability density current which satisfy a
continuity equation.

The derivatives which act on a wave function can be generalized to
derivatives which are covariant under the action of a gauge
transformation with parameters which depend on the lattice point. It
turns out that the corresponding connection can be defined in terms of
a ``moving frame'' which we call, in Section~9, an Einbein. In
Section~10 we discuss the Leibniz rules for these covariant
derivatives. In Section~11 we introduce also a covariant time
derivative and construct a field strength from the commutator of the
two. This permits us to give a differential geometric treatment of our
lattice structure.  In Section~12 it is shown that using our concept
of integration it is possible to derive the Schr\"odinger equation
from a Lagrangian and a variational principle as well as a Noether
theorem.

Finally, in Section~13, we treat the $q$-deformed harmonic oscillator
as a quantum mechanical example in this formalism.
Formally this is similar to the usual treatment of the harmonic
oscillator with creation and annihilation operators~\cite{Lor96}. The wave
functions can be found by solving $q$-differential equations. The
Gauss function on the lattice plays an important role. We define a
Fourier transformation with the $q$-cosine and $q$-sine functions. The
ground state of the harmonic oscillator is a Gaussian in the momentum
space. Under a Fourier transformation it is no longer a Gaussian; it
is the $q$-exponential function. This provides a nice example where a
$q$- Fourier transformation can be carried out explicitly.

\sect{The algebra and its representation}

The $q$-deformed Heisenberg algebra on which our calculus will be based
is a formal $*$-algebra generated by elements $(x,p)$ and an extra
generator $\Lambda$ (the dilatator) which satisfy the commutation
relations
\ba
  \label{i1}
  && q^{\frac 12}xp - q^{-\frac 12}px = i\Lambda,\\
  && \Lambda p = qp\Lambda, \quad \Lambda x = q^{-1}x\Lambda
\ea
Here $q$ is a real number greater than one. The elements $x$ and $p$ are 
assumed to be hermitian and $\Lambda$ to be unitary:
\be
  \label{i2}
  \bar x = x,\quad \bar p = p,\quad \bar \Lambda = \Lambda^{-1}
\ee
This ``bar'' operation is meant to be an algebraic involution and
coincides with complex conjugation on numbers $(\bar{q}=q)$.  This
algebra and its representations have been studied by Hebecker
\etal~\cite{Heb94} and Schm\"udgen~\cite{Sch98}. In these
representations the bar operation is realized as the star operation on
linear operators. We are interested only in those representations
where a formally hermitian operator is represented by a self-adjoint
linear operator on $L_q^2$.  In such representations the operator $x$
can be assumed to be diagonal and its eigenvalues are~\cite{Heb94} given
by
\ba
  \label{i3}
  && x|n,\sigma\ra^s = \sigma sq^n|n,\sigma\ra^s\\
  && n\in\b{Z},\quad\sigma=\pm1\nonumber
\ea
The number $s$ characterizes the representation and can take the
values $1 \leq s < q$. The eigenvectors $|n,\sigma\ra^s$ form an
orthonormal basis of $L^2_q$:
\be
  \label{i4}
  {}^s\la n',\sigma'|n,\sigma\ra^s = \delta_{n,n'}\delta_{\sigma,\sigma'}
\ee
The action of the operator $\Lambda$ in the above representation is given by
\be
  \label{i5}
  \Lambda|n,\sigma\ra^s = |n+1,\sigma\ra^s
\ee

The action of $p$ can now be obtained from (\ref{i3}) and (\ref{i5}).
We must first  enlarge the algebra by adding the element $x^{-1}$. 
This element is well defined on the basis $|n,\sigma\ra^s$.
Next we conjugate the relation (\ref{i1}):
\be
  \label{i6}
  q^{\frac 12}px - q^{-\frac 12}xp = -i\Lambda^{-1}
\ee
When we eliminate $px$ from the equations (\ref{i1}) and (\ref{i6})
we obtain
\be
  \label{i7}
  xp = i\frac{q^{\frac 12}}{(q-q^{-1})}(\Lambda - q^{-1}\Lambda^{-1})
\ee
We multiply this relation by $x^{-1}$ and we find that
\be
  \label{i8}
  p = i\frac{q^{\frac 12}}{(q-q^{-1})}x^{-1}(\Lambda - q^{-1}\Lambda^{-1})
\ee
For the representations (\ref{i3}) and (\ref{i5}) this yields (s=1)
\be
  \label{i9}
  p|n,\sigma\ra = i\frac{\sigma}{(q-q^{-1})}q^{-n}
  \Big(q^{-\frac 12}|n+1,\sigma\ra - q^{\frac 12}|n-1,\sigma\ra\Big)
\ee

We note that $\sigma$ does not change under the action of $x$, $p$ and
$\Lambda$. Therefore each sign of $\sigma$ yields a
representation. But both signs of $\sigma$ must be used to find
self-adjoint extensions of the hermitian operator $p$, as it is
defined by (\ref{i9}), which then also satisfy the relations (\ref{i1}).

The algebraic relations allow an arbitrary ordering of the elements
($x$, $p$, $\Lambda$); any product of these elements, in arbitrary
order, can be expressed in terms of ordered polynomials, for example
in the order $xp\Lambda$. We notice also that the algebra defined by
the relations (\ref{i1}) has a subalgebra which is generated by the
elements ($p, p^{-1}, \Lambda,\Lambda^{-1})$ as well as a subalgebra
generated by the elements ($x, x^{-1}$). Elements of this latter
algebra we shall call fields. By $f(x)$ is meant an element of the
algebra generated by ($x, x^{-1}$) which then is completed by
allowing formal power series.

\sect{Derivatives}

In the previous section we defined fields. Derivatives will be 
mappings of this algebra into itself which we shall now define. 
{From} the ordering property of the
algebra we know that for any field $f$ there are fields $g$ and $h$
such that
\be
  \label{i10}
  pf(x)=g(x)p - iq^{\frac 12}h(x)\Lambda
\ee
as well as a field $j(x)$ such that
\be
  \label{i11}
  \Lambda f(x) = j(x)\Lambda
\ee
A derivative is now defined as the map
\be
  \label{i12}
  \nabla f(x) = h(x)
\ee
In addition we define
\be
  \label{i13}
  Lf(x) = j(x) 
\ee
On the monomials $x^n, (n\in\b{Z})$ these maps are given by
\ba
  \label{i14}
  && \nabla x^n = [n]x^{n-1},\qquad [n] = \frac{q^n-q^{-n}}{q-q^{-1}}\\
  && Lx^n = q^{-n} x^n  \nonumber
\ea
These monomials form a basis of the algebra generated by $(x, x^{-1})$. 
It follows therefore from (\ref{i14}) that there is an
algebra morphism of the $(p, \Lambda)$ algebra to the $(\nabla, L)$
algebra define by the relation
\be
  \label{i15}
  L\nabla = q\nabla L
\ee
It follows also from (\ref{i14}) that the action of $\nabla$ can be
generated by $L, L^{-1}$ and $x^{-1}$:
\be
  \label{i16}
  \nabla = \frac{1}{q-q^{-1}} x^{-1} (L^{-1}-L)
\ee
In this formula $x^{-1}$ is to be interpreted as a map defined by left
multiplication within the $(x, x^{-1})$ algebra.

Next we consider the Leibniz rule. From the formula
\be
  \label{i17}
  x^{m+n} = x^mx^n
\ee
it follows that
\ba
  \label{i18}
  && (Lx^{m+n}) = (Lx^m)(Lx^n)\\
  && (L^{-1}x^{m+n}) = (L^{-1}x^m)(L^{-1}x^n)\nonumber
\ea
{From} the expression (\ref{i16}) for $\nabla$ we obtain from (\ref{i18})
the Leibniz rule
\ba
  \label{i19}
  \nabla(fg) &=& (\nabla f)(Lg) + (L^{-1}f)(\nabla g)\\
             &=& (\nabla f)(L^{-1}g) + (L f)(\nabla g)\nonumber
\ea
for the derivative.  Equations (\ref{i18}) and (\ref{i19})
can be seen as comultiplication rules for the $(\nabla, L)$ algebra:
\ba
  \label{i20}
  \Delta(\nabla) &=& \nabla \otimes L + L^{-1} \otimes \nabla\\
  \Delta(L) &=& L \otimes L\nonumber
\ea
or
\ba
  \label{i21}
  \Delta(\nabla) &=& \nabla \otimes L^{-1} + L \otimes \nabla\\
  \Delta(L) &=& L \otimes L\nonumber
\ea
It is easy to see that this is an algebra morphism. Acting on fields
the two comultiplication rules of $\nabla$ coincide.
The map defined by (\ref{i12}) is not onto since $x^{-1}$ is not in
the image of $\nabla$. The kernel of this map consists of the constants:
\be
  \label{i22}
  \nabla c = 0
\ee

\sect{Differentials}

We define the differential or exterior derivative on elements of the
$(x, x^{-1})$ algebra as
\be
  \label{i23}
  d = dx\,\nabla, \quad dx = (dx)
\ee
It has a unique extension to 1-forms if one imposes $d^2 = 0$. Because
$x$ and $dx$ can be ordered there are no higher-order forms. If we apply
$d$ to $x^n$ we obtain
\be
  \label{i24}
  dx^n = dx\,\nabla x^n = [n]dx\, x^{n-1}
\ee
To find a Leibniz rule we start from the relation:
\be
  \label{i25}
  dx^{m+n} = [m+n]dx\, x^{m+n-1}
\ee
and we try the Ansatz
\be
  \label{i26}
  dx^mx^n = dx^m(Ax^n)+(Bx^m)dx^n
\ee
This can be compared with (\ref{i25}):
\be
  \label{i27}
  [m]dx\, x^{m-1}(Ax^n)+(Bx^m)[n]dx\, x^{n-1} = [m+n]dx\, x^{m+n-1} 
\ee
We collect the terms with $dx$ on the left hand side:
\be
  \label{i28}
  dx\Big([m+n]x^{m+n-1} - [m]x^{m-1}(Ax^n)\Big) = [n] (Bx^m)dx\, x^{n-1}
\ee
We find that the left hand side has to be proportional to $[n]$. This
can be achieved by setting $A=L^{-1}$ or $A=L$. We first analyze the
case $A=L^{-1}$
\be
  \label{i29}
  (Ax^n)=(L^{-1}x^n) = q^n x^n
\ee
Eqn (\ref{i28}) becomes
\be
  \label{i30}
  dx\, q^{-m} x^{m+n-1} = (Bx)^m dx\, x^{n-1}
\ee
or
\be
  \label{i31}
  dx\, x^m = q^m(Bx)^m dx
\ee
If we set $B=L^b$ with $b\in\b{Z}$ we can derive from (\ref{i31}) 
the commutation relation
\be
  \label{i32}
  dx\, x = q^{1-b} x dx
\ee
between $dx$ and $x$. This leads us to the Leibniz rule
\ba
  \label{i33}
  && d(fg) = df(L^{-1}g) + (L^bf)dg \\
  && dx\, x = q^{1-b} x dx \nonumber
\ea
The second choice for $A$ is $A=L$. The equation analogous to (\ref{i30}) is 
now
\be
  \label{i34}
  dx\, q^m x^{m+n-1} = (Bx)^m dx\, x^{n-1}
\ee
or
\be
  \label{i35}
  dx\, x^m = q^{-m} (Bx)^m dx.
\ee
We obtain the Leibniz rule:
\ba
  \label{i36}
  && d(fg) = df(Lg) + (L^bf)dg \\
  && dx\, x = q^{-1-b} x dx \nonumber
\ea
An interesting choice for the Leibniz rule (\ref{i33}) is $b=1$ or for
(\ref{i36}) $b=-1$; in both cases $dx$ and $x$ commute.

\sect{The integral}

We define the indefinite integral to be the inverse image of
(\ref{i12}).  Integrals have also been defined by A. Kempf and S. Majid in 
\cite{Kem94} and by H. Steinacker in \cite{Ste96}.
The kernel of the map (\ref{i12}) are the constants and
$x^{-1}$ is not in the image of $\nabla$. Thus we find
\be
  \label{i37}
  \int^x x^n = \frac{1}{[n+1]} x^{n+1} + \mbox{const}.\qquad 
  n\in\b{Z},\, n\neq -1
\ee
A useful formula is obtained if we invert $\nabla$ in the form
(\ref{i16}):
\be
  \label{i38}
  \nabla^{-1} = (q-q^{-1}) \frac{1}{L^{-1}-L}\, x
\ee
For $m\neq -1$ we can apply this to $x^m$ and obtain
\ba
  \label{i39}
  \nabla^{-1}x^m &=& (q-q^{-1}) \frac{1}{q^{m+1}-q^{-m-1}}\, x^{m+1}\\
                 &=& \frac{1}{[m+1]}\, x^{m+1} \nonumber
\ea
If we apply (\ref{i38}) to a field we can expand:
\ba
  \label{i40}
  \nabla^{-1}f(x) &=& (q-q^{-1})\sum_{\nu = 0}^{\infty} 
                      L^{2\nu}Lxf(x)\\
                  &=& -(q-q^{-1})\sum_{\nu = 0}^{\infty} 
                      L^{-2\nu}L^{-1}xf(x)\nonumber
\ea
These two formulas should be used depending what series converges.
That is, for $x^m, m\ge 0$ the first series converges, for $x^m, m<-1$ the
second series converges.

A definite integral can be defined only once a representation of the
algebra (\ref{i1}) is given. We consider a representation where the
$s$ of Equation~(\ref{i3}) is equal to one $(s=1)$. The linear
operator $x$ has eigenvalues $\sigma q^M$. Let us first consider the
case $\sigma = +1$. A definite integral would be an integral from
$q^N$ to $q^M$. It should be in agreement with (\ref{i37}) for
monomials:
\be
  \label{i41}
  \int\limits_N^M x^n = \frac{1}{[n+1]}\Big(q^{M(n+1)}-q^{N(n+1)}\Big)
\ee
For a general field we can extend it by linearity. This definition has 
Stokes' theorem as a consequence:
\be
  \label{i42}
  \int\limits_N^M \nabla x^n = q^{Mn}-q^{Nn} = x^n\Big|_N^M
\ee
The definition (\ref{i41}) is not suitable to define an integral over
a function in the limit $N\to -\infty, M\to\infty$.  To define such an
integral we start from Equation~(\ref{i40}) and apply it to $x^m,
m>0$. We see that the powers of $L$ take even values. Therefore we
shall study the integral (\ref{i42}) with even and odd $M$ separately,
taking immediately the limit $N\to -\infty$.
\ba
  \label{i43}
  \int\limits_{-\infty}^{2M} x^m &=& (q-q^{-1}) 
  \sum_{\nu=0}^{\infty} q^{-(2\nu+1)(m+1)}q^{2M(m+1)}\\
  &=& \frac{1}{[m+1]}q^{2M(m+1)}\nonumber 
\ea
This agrees with (\ref{i41}). Next we rewrite the sum in (\ref{i43})
\ba
  \label{i44}
  \int\limits_{-\infty}^{2M} x^m &=& 
  (q-q^{-1})\sum_{\nu=0}^{\infty} q^{(m+1)(2M-2\nu-1)}\nonumber\\
  &=& (q-q^{-1})\sum_{\mu=-\infty}^{M} q^{(m+1)(2\mu-1)}\\
  &=& (q-q^{-1})\sum_{\mu=-\infty}^{M} \la 2\mu|Lxx^m|2\mu\ra\nonumber
\ea
where $|2\mu\ra$ are states of the representation (\ref{i3}). We use
this formula for polynomials $h(x)$:
\be
  \label{i45}
  \int\limits_{-\infty}^{2M}h(x) = 
  (q-q^{-1})\sum_{\mu=-\infty}^{M} \la 2\mu|Lxh(x)|2\mu\ra
\ee
In a similar way we find
\be
  \label{i46}
  \int\limits_{-\infty}^{2M+1}h(x) = 
  (q-q^{-1}) \sum_{\mu=-\infty}^{M} \la 2\mu+1|Lxh(x)|2\mu+1\ra
\ee
For negative powers of $x (x^m, m\le -2)$ we use the second sum in
(\ref{i40}) and find for the respective polynomials:
\ba
  \label{i47}
  \int\limits_{2M}^{\infty}h(x) &=& 
  (q-q^{-1})\sum_{\mu=M+1}^{\infty} \la 2\mu|Lxh(x)|2\mu\ra\\
  \int\limits_{2M+1}^{\infty}h(x) &=& 
  (q-q^{-1})\sum_{\mu=M+1}^{\infty} \la 2\mu+1|Lxh(x)|2\mu+1\ra\nonumber
\ea
It is now obvious how a definite integral for a field should be formulated:
\ba
 \label{i48}
  \int\limits_{2N}^{2M}h(x) &=& 
  (q-q^{-1})\sum_{\mu=N+1}^{M} \la 2\mu|Lxh(x)|2\mu\ra\\
  \int\limits_{2N+1}^{2M+1}h(x) &=& 
  (q-q^{-1})\sum_{\mu=N+1}^{M} \la 2\mu+1|Lxh(x)|2\mu+1\ra\nonumber
\ea
It is part of an interesting structure that also manifests itself in
the $q$-Fourier transform that odd and even valued lattice points are
quite independent. Thus we have two representations of the definite
integral over odd or even points.  The part of the representation
(\ref{i3}) with $\sigma=-1$ can be treated completely analogously, but
again even and odd are quite independent.
For monomials (\ref{i48}) is identical with (\ref{i41}). This shows
that Stokes' theorem holds for even and odd $M,N$ as well:
\be
  \label{i49}
  \int\limits_N^M \nabla f(x) = f(x)\Big|_{q^N}^{q^M}
\ee
If for both integrals (\ref{i48}) the limit 
$M\to \infty, N\to -\infty$ exists as well as the corresponding 
integral over negative eigenvalues of $x$, then we define the integral as
\ba
  \label{i50}
  \int h(x) &=& \frac 12(q-q^{-1})\sum_{\sigma=+,-}
  \lim_{M\to\infty \atop N\to -\infty}
  \sigma\,\Big\{\sum_{\mu=N}^M\la2\mu,\sigma|Lxh(x)|2\mu,\sigma\ra\nonumber\\
  &&+\la 2\mu+1,\sigma|Lxh(x)|2\mu+1,\sigma\ra\Big\}\\
  &=& \frac 12(q-q^{-1})\sum_{\sigma=+,-}\sigma
  \sum_{\mu=-\infty}^{\infty} \la\mu,\sigma|Lxh(x)|\mu,\sigma\ra\nonumber
\ea
We assume that the rearrangement of the sum is allowed.  

It remains to define the integral over $x^{-1}$. For the definite integral
(\ref{i48}) this is possible:
\be
  \label{i51}
  \int\limits_{2N}^{2M}\frac{1}{x} = (q-q^{-1})(M-N)
\ee
For the integration limits  $\bar{z}=q^{2M}, \underline{z}=q^{2N}$ formula
(\ref{i51}) approaches $\ln\bar{z}-\ln\underline{z}$ for $q\to 1$.

If $h(x)\to 0$ for $x\to\pm\infty$ we conclude from
 (\ref{i49})
\be
  \label{i52}
  \int\nabla h(x) =0
\ee

\sect{The Hilbert space $L^2_q$}

The integral (\ref{i50}) can be used to define a scalar product:
\be
  \label{i53}
  (\chi,\psi) = \int \chi^{\ast}\psi = 
  \frac{q-q^{-1}}{2}\sum_{\sigma=+1,-1}
  \sum_{\mu=-\infty}^{\infty} 
  \sigma \la\mu,\sigma|Lx\chi^{\ast}\psi|\mu,\sigma\ra
\ee
The factor $Lx$ in the matrix element allows the sum to converge at
$x=0 (q^N,N\rightarrow-\infty)$ for fields that do not vanish at
$x=0$. It is the same factor that occurs in the Jackson integral. It
should be noted, however, that the integral (\ref{i53}) has been
obtained in Equation~(\ref{i50}) from a sum over even and odd values of
$\mu$ seperately.

For reasonably behaved fields we conclude from Stokes' theorem (\ref{i52}) 
that 
\be
  \label{i54} \int \nabla(\chi^{\ast}\psi) = 0
\ee
and find from (\ref{i19}) we find a formula for partial integration:
\ba
  \label{i55}
  \int (\nabla\chi^{\ast})(L\psi) + 
  \int (L^{-1}\chi^{\ast})(\nabla\psi) = 0\\
  \int (\nabla\chi^{\ast})(L^{-1}\psi) + 
  \int (L\chi^{\ast})(\nabla\psi) = 0\nonumber
\ea
We use these formulas to find a $q$-version of Green's theorem. We write:
\ba
  \label{i56}
  \nabla\big((\nabla f)(L^{-1}g)\big) &=& 
  (\nabla^2 f)(g) + (L^{-1}\nabla f)(\nabla L^{-1}g)\\
  \nabla\big((L^{-1} f)(\nabla g)\big) &=& 
  (\nabla L^{-1} f)(L^{-1}\nabla g) + (f)(\nabla^2 g)\nonumber
\ea
We subtract these two equations and obtain Green's theorem:
\be
  \label{i57}
  (\nabla^2 f)(g) - (f)(\nabla^2 g) = 
  \nabla\Big((\nabla f)(L^{-1}g) - (L^{-1}f)(\nabla g)\Big)
\ee
As a consequence of (\ref{i54}) we find that $\nabla^2$ is an hermitean
operator:
\be
  \label{i58}
  \int (\nabla^2\chi^{\ast})\psi = \int (\chi^{\ast})(\nabla^2\psi)
\ee
If we integrate (\ref{i57}) over a finite volume we find
\be
  \label{i59}
  \int\limits_N^M \big\{(\nabla^2f)(g)-(f)(\nabla^2g)\big\} = 
  \big\{(\nabla f)(L^{-1}g)-(L^{-1}f)(\nabla g)\big\}\Big|_N^M
\ee
To define $\nabla^2$ we have implicitely used a metric. For a further 
discussion of this point we refer to Cerchiai \etal~\cite{Cer98}.

\sect{The $\cos_q$ and $\sin_q$ functions:}

The $q$-deformed cosine and sine functions are defined as follows:
\ba
  \label{i60}
  \sin_q(z) &=& \sum_{n=0}^{\infty} 
  (-1)^n q^{-2n(n+1)}\frac{z^{2n+1}}{(q^{-2};q^{-2})_{2n+1}}\\
  \cos_q(z) &=& \sum_{n=0}^{\infty} 
  (-1)^n q^{-2n(n+1)}\frac{z^{2n}}{(q^{-2};q^{-2})_{2n}}\nonumber
\ea
with:
\be
  \label{i61}
  (q^{-2};q^{-2})_{n} = q^{-\frac 12n(n+1)}\big(q-q^{-1}\big)^n[n]!
\ee
These functions are well behaved at the points $z=q^{2l}$ and satisfy
an orthogonality and completeness relation:
\ba
  \label{i62}
  \sum_{k=-\infty}^{\infty} q^{-2k}\sin_q(q^{-2(k+n)})\sin_q(q^{-2(k+m)}) = 
   \frac{q^{2m}}{N_q^2}\delta_{n,m}\\
  \sum_{k=-\infty}^{\infty} q^{-2k}\cos_q(q^{-2(k+n)})\cos_q(q^{-2(k+m)}) = 
  \frac{q^{2m}}{N_q^2}\delta_{n,m}\nonumber
\ea
with 
$$
N_q = \frac{(q^{-2};q^{-4})_{\infty}}{(q^{-4};q^{-4})_{\infty}}.
$$
This formula shows that $\cos_q(q^{2l})$ and $\sin_q(q^{2l})$ have to
vanish strong enough for $l\to +\infty$ to make the sum in (\ref{i62})
converge. This would not be the case for odd powers of $q$.

Formulas (\ref{i62}) tell us that a function $f$ can be expanded in
terms of the $\cos_q$ $\underline{\hbox{or}}$ the $\sin_q$ functions.
\ba
  \label{i63}
  g(q^{-2n}) &=& N_q\sum_{k=-\infty}^{\infty}
  q^{-2k}\cos_q(q^{-2(k+n)})f(q^{-2k})\\ 
  g(q^{-2n+1}) &=& N_q\sum_{k=-\infty}^{\infty}
  q^{-2k}\cos_q(q^{-2(k+n)})f(q^{-2k+1})\nonumber 
\ea
The functions $g$ are the expansion coefficients of $f$ and $f$ can be
obtained from $g$ with the help of the orthogonality relation
(\ref{i62})
\ba
  \label{i64} 
  f(q^{-2k}) &=& N_q\sum_{n=-\infty}^{\infty}
  q^{-2n}\cos_q(q^{-2(k+n)})g(q^{-2n})\\ 
  f(q^{-2k+1}) &=& N_q\sum_{n=-\infty}^{\infty}
  q^{-2n}\cos_q(q^{-2(k+n)})g(q^{-2n+1})\nonumber 
\ea
We see again that even and odd powers of $q$ play a very independent role.
{From} (\ref{i62}) follows that the transformation (\ref{i63}) is an 
isometry:
\ba
  \label{i65}
  \sum_k q^{-2k}f^{\ast}(q^{-2k})f(q^{-2k}) &=& 
  \sum_n q^{-2n}g^{\ast}(q^{-2n})g(q^{-2n})\\
  \sum_k q^{-2k+1}f^{\ast}(q^{-2k+1})f(q^{-2k+1}) &=& 
  \sum_n q^{-2n+1}g^{\ast}(q^{-2n+1})g(q^{-2n+1})\nonumber
\ea
A similar analysis could have been carried out with the $\sin_q$ function.

There is an important relation between $\cos_q$ and $\sin_q$:
\ba
  \label{i66}
  \frac{1}{z}\big(\cos_q(z)-\cos_q(q^{-2}z)\big) &=& -q^{-2}\sin_q(q^{-2}z)\\
  \frac{1}{z}\big(\sin_q(z)-\sin_q(q^{-2}z)\big) &=& \cos_q(z)\nonumber
\ea
If we compare this formula with our formula (\ref{i16}) for $\nabla$
we find that
\ba
  \label{i67}
  \nabla\cos_q(x) &=& -\frac{1}{q}\frac{1}{(q-q^{-1})}L\sin_q(x)\\
  \nabla\sin_q(x) &=& \frac{q}{(q-q^{-1})}L^{-1}\cos_q(x)\nonumber
\ea

Now we want to construct a basis in $L^2_q$.
We split $L^2_q$ into four subspaces with $\sigma=+1$, $\sigma =-1$
and $\mu$ even or $\mu$ odd.
\be
  \label{i68}
  \c{H} = \c{H}^{\mu\,\mbox{\scriptsize even}}_{\sigma =
   +1}\oplus\c{H}^{\mu\,\mbox{\scriptsize odd}}_{\sigma=
   +1}\oplus\c{H}^{\mu\,\mbox{\scriptsize even}}_{\sigma=
   -1}\oplus\c{H}^{\mu\,\mbox{\scriptsize odd}}_{\sigma=-1}
\ee
On each of these subspaces we define a projector
$\Pi^{\mbox{\scriptsize even}}_+$, 
$\Pi^{\mbox{\scriptsize odd}}_+$, 
$\Pi^{\mbox{\scriptsize even}}_-$ and $\Pi^{\mbox{\scriptsize odd}}_-$
respectively. We start by defining a basis in
$\c{H}^{\mu\,\mbox{\scriptsize even}}_{\sigma=+1}$. The norm in this
subspace is:
\be
  \label{i69}
  (\psi,\psi) = \frac 12(q-q^{-1})\sum_{\mu-\infty}^{\infty}
  \la 2\mu,\sigma=+1|Lx\psi^{\ast}\psi|2\mu,\sigma=+1\ra
\ee
The functions:
\be
  \label{i70}
  \c{C}^{(+)}_{2n+1 \atop n\in\b{Z}}(x) = 
  \tilde{N}_q \cos_q(xq^{2n+1})\Pi^{\mbox{\scriptsize even}}_+,\qquad 
  \tilde{N}_q = N_q\Big(\frac{2q}{q-q^{-1}}\Big)^{\frac 12}q^n 
\ee
form an orthonormal basis in 
$\c{H}^{\mu\,\mbox{\scriptsize even}}_{\sigma=+1}$.  
This follows from (\ref{i63}) and (\ref{i64}) and the definition of
the scalar product (\ref{i69}). Similarly we define
\be
  \label{i71}
  \c{C}^{(+)}_{2n}(x) = 
  \tilde{N}_q \cos_q(xq^{2n})\Pi^{\mbox{\scriptsize odd}}_+
\ee
in the subspace $\c{H}^{\mu\,\mbox{\scriptsize odd}}_{\sigma=+1}$ with
the norm
\be
  \label{i72}
  (\psi,\psi) = 
  \frac 12(q-q^{-1})\sum_{\mu-\infty}^{\infty}
  \la 2\mu+1,\sigma=+1|Lx\psi^{\ast}\psi|2\mu+1,\sigma=+1\ra
\ee
and 
\ba
  \label{i73}
  \c{C}^{(-)}_{2n+1}(x) &=& 
  \tilde{N}_q\cos_q(xq^{2n+1})\Pi^{\mbox{\scriptsize even}}_-\\
  \c{C}^{(-)}_{2n}(x) &=& 
  \tilde{N}_q\cos_q(xq^{2n})\Pi^{\mbox{\scriptsize odd}}_-\nonumber
\ea
In an analogous way we could have used $\sin_q$ in any of the
subspaces to define $\c{S}^{(+)}_{2n+1}(x)$, $\c{S}^{(+)}_{2n}(x)$,
$\c{S}^{(-)}_{2n+1}(x)$ and $\c{S}^{(-)}_{2n}(x)$ any combination of
the individual basis can be used to define a basis in all of $L_q^2$.
We note, however, that the operator $\nabla$ is not defined on any
of the elements of this basis. It maps elements of the basis into 
functions that are not in $L^2_q$.  In contrast $\nabla L$ or $\nabla
L^{-1}$ is defined on this basis. These are exactly the operators that
enter in Green's theorem (\ref{i57}) if we write it in the form:
\be
  \label{i74}
  (\nabla^2f)(g) - (f)(\nabla^2 g) = 
  \nabla L^{-1}\big\{(L\nabla f)(g)-(f)(L\nabla g)\big\}
\ee

The elements of our basis are eigenfunctions of the operator
$\nabla^2$. To see this we generalize (\ref{i67})
\ba
  \label{i75}
  \nabla_{\!x}\cos_q(xy) &=& 
  \frac{1}{q-q^{-1}}x^{-1}\big(\cos_q(qxy)-\cos_q(q^{-1}xy)\big)\nonumber\\
                     &=& -\frac{1}{q-q^{-1}} q^{-1} y L \sin_q(xy)\\
  \nabla_{\!x}\sin_q(xy) &=& 
  \frac{1}{q-q^{-1}} q y L^{-1}\cos_q(xy)\nonumber
\ea
and we find
\ba
  \label{i76}
  \nabla^2\cos_q(xy) &=& 
  -\frac{1}{q}\big(q-q^{-1}\big)^{-2} y^2 \cos_q(xy)\\
  \nabla^2\sin_q(xy) &=& 
  -q\big(q-q^{-1}\big)^{-2} y^2 \sin_q(xy)\nonumber
\ea
We see that the basis represents eigenfunctions with the following
eigenvalues:
\ba
  \label{i77}
  \c{C}^{+,-}_{2n+1} &:& -\big(q-q^{-1}\big)^{-2}q^{4n+1}\nonumber\\
  \c{C}^{+,-}_{2n}   &:& -\big(q-q^{-1}\big)^{-2}q^{4n-1}\\
  \c{S}^{+,-}_{2n+1} &:& -\big(q-q^{-1}\big)^{-2}q^{4n+3}\nonumber\\
  \c{S}^{+,-}_{2n}   &:& -\big(q-q^{-1}\big)^{-2}q^{4n+1}\nonumber
\ea
We also see that the set of eigenfunctions is overcomplete. We
conclude that $\nabla^2$ is hermitian but not self-adjoint. On any of
the bases defined before a self-adjoint extension of the operator
$\nabla^2$ is defined. Note that the set of eigenvalues does depend on
the respective extension.

\sect{Schr\"odinger equation}

It is natural to define the Schr\"odinger equation as follows:
\be
  \label{ii1}
  i\frac{\partial}{\partial t}\psi = \big(-\frac{1}{2m}\nabla^2+V\big)\psi
\ee
It has a probability interpretation. To show this we calculate:
\be
  \label{ii2}
  \frac{d}{dt}\psi^{\ast}\psi = - \frac{1}{2mi}\big\{\psi^{\ast}(\nabla^2\psi)-
 (\nabla^2 \psi^{\ast})\psi\big\}
\ee
We can use (\ref{i57}):
\ba
  \label{ii3}
  \frac{d}{dt}\psi^{\ast}\psi &=& 
   -\frac{1}{2mi}\nabla\big\{(L^{-1}\psi^{\ast})
   (\nabla \psi)-(\nabla\psi^{\ast})(L^{-1}\psi)\big\}\\
   &=& -\frac{1}{2mi}\nabla L^{-1}\big\{\psi^{\ast}
   (L\nabla\psi)-(L\nabla\psi^{\ast})\psi\big\}\nonumber
\ea
This is the continuity equation for the density
\ba
  \label{ii4}
  &&\rho = \psi^{\ast}\psi\\
  &&j = \frac{1}{2mi}L^{-1}\big\{\psi^{\ast}
  (L\nabla\psi)-(L\nabla\psi^{\ast})\psi\big\}\nonumber\\
  &&\frac{\partial}{\partial t}\rho + \nabla j =0\nonumber
\ea
If we integrate this equation over a finite volume we obtain from
(\ref{i49}):
\be
  \label{ii5}
  \frac{d}{dt}\int\limits_{2N}^{2M}\psi^{\ast}\psi = 
   - \frac{1}{2mi}\big\{L^{-1}\big(\psi^{\ast}
   (L\nabla\psi)-(L\nabla\psi^{\ast})\psi\big)\big\}\Big|_{2N}^{2M}
\ee
Note that the ``velocity'' operator is proportional to $L\nabla$.  If
we consider the free Schr\"odinger equation we know the solutions of
the time-indepen\-dent equation:
\be
  \label{ii6}
  -\frac{1}{2m} \nabla^2\psi = E\psi
\ee
These are the functions $\c{S}$ and $\c{C}$ defined in
eqn. (\ref{i70}), (\ref{i71}) and (\ref{i73}). The eigenvalues can be
found in (\ref{i77}).

As an example let us consider the eigenfunction $\c{C}^{(+)}_{2n+1}$.
The energy eigenvalue is
\be
  \label{ii7}
  E^{(+)}_{2n+1} = \frac{(2m)^{-1}}{\big(q-q^{-1}\big)^2}q^{4n+1}
\ee
The normalized eigenfunction is 
\be
  \label{ii8}
  \psi^{(+)}_{2n+1} = 
  \tilde{N}_q\cos_q(xq^{2n+1})\Pi^{\mbox{\scriptsize even}}_+
\ee
For the probability of finding the ``particle'' in the volume between
$q^{2N}$ and $q^{2M}$ we obtain
\ba
  \label{ii9}
  && \int\limits_{2N}^{2M}\tilde{N}_q^2\cos_q^2(xq^{2n+1})\nonumber\\
  && = 2N_q^2q^{2n+1}\sum_{\mu=N+1}^{M}
  \la 2\mu|Lx\cos_q^2(xq^{2n+1})|2\mu\ra\\
  && = 2N_q^2q^{2n}\sum_{\mu=N+1}^{M} q^{2\mu}\cos_q^2(q^{2\mu+2n})\nonumber
\ea
We see that the $\cos_q$ function is well defined for even powers of
$q$.  For the time derivative we evaluate the right hand side of
Equation~(\ref{ii2}).  A typical term is:
\be
  \label{ii10}
  L^{-1}\big(\psi^{\ast}(L\nabla\psi)-
  (L\nabla\psi^{\ast})\psi\big)\Big|_{\mbox{\scriptsize at}\, x = q^{2M}}
\ee
We use (\ref{i67}) and again find that the $\cos_q$ and $\sin_q$
functions are well defined and that the flux is zero. The
eigenfunctions (\ref{ii8}) are stationary.

\sect{Covariant derivatives}

We assume that the field $\psi(x)$ transforms under a gauge
transformation as follows:
\be
  \label{ii11}
  \psi'(x) = e^{i\alpha(x)}\psi(x)
\ee
This is not necessarily an abelian gauge group.  We would like to
define a derivative $\c{D}\psi$, such that
\be
  \label{ii12}
  (\c{D}\psi)' = e^{i\alpha}\c{D}\psi
\ee
For an ordinary derivative we know from the generalized Leibniz rule
(\ref{i19}) that: 
\be
  \label{ii13}
  \nabla\psi' = \big(\nabla e^{i\alpha}\big)(L\psi) + 
  \big(L^{-1}e^{i\alpha}\big)\nabla\psi
\ee
For a covariant derivative we make the Ansatz:
\be
  \label{ii14}
  \c{D}\psi = E(\nabla + \phi)\psi
\ee
The field $E$ plays the role of an Einbein and $\phi$ the role of a
connection.  We now determine the transformation law of $E$ and $\phi$
from (\ref{ii11}) and (\ref{ii12}):
\be
  \label{ii15}
  E'(\nabla +\phi')e^{i\alpha}\psi = e^{i\alpha}E(\nabla + \phi)\psi
\ee
We find
\be
  \label{ii16}
  E' = e^{i\alpha}E\big(L^{-1}e^{-i\alpha}\big)
\ee
and
\be
  \label{ii17}
  \phi' = \big(L^{-1}e^{i\alpha}\big)\phi\big(e^{-i\alpha}\big) - 
  \big(\nabla e^{i\alpha}\big)L e^{-i\alpha}
\ee
Note that the last term is without bracket, the ``operator'' $L$ acts
on the field that multiplies $\phi'$ as well. Another way of saying
this is that $\phi$ is $L$-valued:
\be
  \label{ii18}
  \phi = \varphi L
\ee
For $\varphi$ we find from (\ref{ii17}):
\be
  \label{ii19}
  \varphi' = \big(L^{-1}e^{i\alpha}\big)\varphi
   \big(Le^{-i\alpha}\big) - 
   \big(\nabla e^{i\alpha}\big)\big(L e^{-i\alpha}\big)
\ee
We would like to show now that $\varphi$ can be considered to be a
function of $E$.  To do this we recall formula (\ref{i16}) and
define a covariant  operator $\c{L}$ by
\be
  \label{ii20}
  (\c{L}\psi)' = e^{i\alpha}\c{L}\psi
\ee
{From} the comultiplication law of $L$ follows that
\be
  \label{ii21}
  L\psi' = \big(Le^{i\alpha}\big)(L\psi)
\ee
We make the Ansatz:
\be
  \label{ii22}
  \c{L}\psi = \tilde{E}L\psi
\ee
and find the transformation law of $\tilde{E}$ from (\ref{ii11}) and
(\ref{ii20}):
\be
  \label{ii23}
  \tilde{E}' = e^{i\alpha}\tilde{E}\big(Le^{-i\alpha}\big)
\ee
Similarly we define:
\be
  \label{ii24}
  \tilde\c{L}\psi = EL^{-1}\psi
\ee
and find
\be
  \label{ii25}
  E' = e^{i\alpha} E\big(L^{-1}e^{-i\alpha}\big)
\ee
This transformation law agrees with (\ref{ii16}) and this justifies
our choice of $E$ in the definition (\ref{ii24}).

The inverse Einbein $E^{-1}$ transforms as follows:
\be
  \label{ii26}
  {E^{-1}}' = \big(L^{-1}e^{i\alpha}\big)E^{-1}e^{-i\alpha}
\ee
We see that it is a consistent assumption to postulate
\be
  \label{ii27}
  \tilde{E} = \big(LE^{-1}\big)
\ee
With this assumption we find that
\be
               \tilde\c{L}   = \c{L}^{-1}
\ee
We generalize formula (\ref{i16}):
\be
  \label{ii28}
  \c{D} = \frac{1}{\lambda}x^{-1}\big(\tilde\c{L} - \c{L}\big) 
\ee
and we know from (\ref{ii23}) and (\ref{ii26}) that this definition of
$\c{D}$ has the right transformation property. We now rewrite
(\ref{ii28}) to compare it with (\ref{ii14})
\ba
  \label{ii29}
  \c{D} &=& \frac{1}{\lambda}x^{-1}\big\{E(L^{-1}-L) + (E-\tilde{E})L\big\}\\
           &=& E\nabla + \frac{1}{\lambda}x^{-1}(E-\tilde{E})L\nonumber
\ea
We again see that the connection is $L$ valued and find:
\be
  \label{ii30}
  \varphi = \frac{1}{\lambda}x^{-1}(1-E^{-1}\tilde{E})
\ee
It takes a small calculation to verify that $\varphi$ has the
transformation property (\ref{ii19}) as a consequence of (\ref{ii23})
and (\ref{ii25}).  It is interesting to note that the covariant
derivative of the Einbein $E$ vanishes if we choose (\ref{ii30}) for
the connection. We start from a field $H$ that transforms like $E$:
\be
  \label{ii31}
  H' = e^{i\alpha}H\big(L^{-1}e^{-i\alpha}\big)
\ee
Let us first compute $\c{L}$ and $\tilde\c{L}$ applied to $H$:
\ba
  \label{ii32}
  \c{L}H &=& \tilde{E}(LH)\big(L^{-1}\tilde{E}^{-1}\big)\\
  \tilde\c{L}H &=& E\big(L^{-1}H\big)\big(L^{-1}E^{-1}\big)\nonumber
\ea
We now identify $\tilde{E}$ with $\big(LE^{-1}\big)$ as in (\ref{ii27}) and find
\be
  \label{ii33}
  \c{L}H = \big(LE^{-1}\big)(LH)(E)
\ee
Now we identify $H$ with $E$ and obtain
\be
  \label{ii34}
  \c{L} E = E,\qquad \tilde\c{L} E = E
\ee
Combining this with a covariant derivative as in (\ref{ii28}) yields
\be
  \label{ii35}
  \c{D} E = 0
\ee

\sect{Leibniz rule for covariant derivatives}

We first multiply the field $\psi$ by a scalar field $f$:
\be
  \label{ii36}
  f'=f,\quad \c{D}f=\nabla f,\quad \c{L}f=Lf,\quad \tilde\c{L}f=L^{-1}f
\ee
and we compute
\be
  \label{ii37}
  \c{L}f\psi=\tilde{E}Lf\psi=\tilde{E}(Lf)(L\psi)=\c{L}f\c{L}\psi
\ee
Thus we find for $\c{L}$ the same comultiplication rule as for $L$.
Next we apply this comultiplication rule to the scalar $\psi^{\ast}\psi$:
\ba
  \label{ii38}
  \c{L} \psi^{\ast}\psi &=& L(\psi^{\ast}\psi)=
  (L\psi^{\ast})(L\psi)\\
  &=& (L\psi^{\ast})\tilde{E}^{-1}(\tilde{E}L\psi)=
  \c{L}\psi^{\ast}\c{L}\psi\nonumber
\ea
and we conclude that
\be
  \label{ii39} \c{L}\psi^{\ast} = 
  (L\psi^{\ast})\tilde{E}^{-1} = (L\psi^{\ast})(LE)
\ee
We have here used (\ref{ii27}).

If we now apply the comultiplication rule to the product of two
arbitrary representations we find:
\ba
  \label{ii40}
  \c{L}\psi\chi &=& \c{L}\psi\c{L}\chi\\
  &=& (\tilde{E}L)\psi(\tilde{E}L\chi)\nonumber
\ea
It becomes more transparent if we use indices:
\ba
  \label{ii41}
  \c{L}\psi_{\mu}\chi_{\alpha} &=& 
  \tilde{E}_{\mu}{}^{\rho}\tilde{E}_{\alpha}{}^{\beta}
  (L\psi)_{\rho}(L\psi)_{\beta}\\
  &=& \tilde{E}_{\mu\alpha}{}^{\rho\beta}
  L\psi_{\rho}\chi_{\beta}\nonumber
\ea
Thus the Einbein matrix for the product representation is the product
of the Einbein matrices of the respective representations.  Therefore
we can generate the Einbein for any representation from the Einbein of
the fundamental representation.  This immediately generalizes to
$\tilde\c{L}$ and $\c{D}$. We find the comultiplication rule
\be
  \label{ii42}
  \tilde\c{L}\psi\chi = \tilde\c{L}\psi\tilde\c{L}\chi  
\ee
and the Leibniz rule
\be
  \label{ii43}
  \c{D}\psi\chi = (\c{D}\psi)\c{L}\chi + (\tilde\c{L}\psi)\c{D}\chi 
\ee
It should be noted that the transformation parameter $\alpha$ in
(\ref{ii11}) is proportional to a coupling constant $g$. Thus we
should expand in the coupling constant:
\be
  \label{ii44}
  E_{\alpha}{}^{\beta} = \delta_{\alpha}{}^{\beta} + g h_{\alpha}{}^{\beta}(x)
\ee

\sect{Curvature}

To be able to speak about curvature we must first define a covariant time
derivative:
\be
  \label{ii45}
  \c{D}_t\psi = (\partial_t +\omega)\psi
\ee
The connection $\omega$ transforms as follows:
\be
  \label{ii46}
  \omega' = e^{i\alpha}\omega e^{-i\alpha}+e^{i\alpha}\partial_te^{-i\alpha}
\ee
Now we can compute the commutator of the two covariant derivatives:
\ba
  \label{ii47}
  (\c{D}_t \c{D} - \c{D} \c{D}_t)\psi &=& 
  \big\{(\partial_t+\omega)E(\nabla+\phi) - 
  E(\nabla+\phi)(\partial_t+\omega)\big\}\psi\nonumber\\
  &=& \big\{(\partial_tE)E^{-1}-E(L^{-1}\omega)E^{-1}+ 
  \omega\big\}\c{D}\psi\\
  && + \big\{\partial_tE\phi-E(\nabla\omega)L-
  E\phi\omega-(\partial_tE)\phi+E(L^{-1}\omega)\phi\big\}\psi\nonumber
\ea
This allows us to define two tensors:
\ba
  \label{ii48}
  T &=& (\partial_tE)E^{-1}-(E)(L^{-1}\omega)E^{-1} + \omega\\
  F &=& \big\{\partial_t\varphi - \nabla\omega + 
  (L^{-1}\omega)\varphi - \varphi(L\omega)\big\}\nonumber
\ea
and to write (\ref{ii47}) in the form
\be
  \label{ii49}
  (\c{D}_t \c{D} - \c{D} \c{D}_t)\psi = T\c{D}\psi + EFL\psi
\ee

Starting from the transformation laws of $E$, $\omega$ and $\varphi$,
these are the eqn. (\ref{ii16}), (\ref{ii19}) and (\ref{ii46}). We can
verify by a lengthy calculation the transformation laws of $T$ and $F$
as they follow from the definition (\ref{ii49}):
\ba
  \label{ii50}
  T' &=& e^{i\alpha}Te^{-i\alpha}\\
  F' &=& \big(L^{-1}e^{i\alpha}\big)F\big(Le^{-i\alpha}\big)\nonumber
\ea
The tensor $T$ has already the right transformation law. It follows from 
(\ref{ii16}) that the quantity
\be
  \label{ii51}
  \c{F} = EF(LE)
\ee
transforms as a tensor as well: 
\be
  \label{ii52}
  \c{F}' = e^{i\alpha}\c{F}e^{-i\alpha}
\ee
The curvature $T$ can also be derived from the following commutation:
\be
                 (\c{L}\c{D}_t - \c{D}_t\c{L}) \psi = LE^{-1}T\psi
\ee

\sect{Euler-Lagrange equation and Noether theorem}

The definition of the integral (\ref{i48}) allows the formulation of a
variational problem. An action can be defined as the integral over a
Lagrangian. The Lagrangian itself is a function of the fields and
their derivatives and we demand that the action be extremal under the
variation of the fields with fixed boundary values. Let us first
examine the Schr\"odinger equation (\ref{ii1}). We define the action
\be
  \label{ii53}
  W=\int\limits_{t_1}^{t_2}dt\, \int\limits_{2N}^{2M}
  \psi^{\ast}\big(i\frac{\partial}{\partial t}\psi + 
  \frac{1}{2m}\nabla^2\psi - V\psi\big)
\ee
Variation of $\psi^{\ast}$ leads to the Schr\"odinger equation
(\ref{ii1}). We know that $\nabla^2$ is hermitian. Thus a variation of
$\psi$ leads to the conjugate Schr\"odinger equation.

We would prefer to formulate the Lagrangian in terms of the fields and
their first derivatives. To write (\ref{ii53}) in such a form we note
that the Leibniz rule (\ref{i19}) implies that
\be
  \label{ii54}
  \int (\nabla L\psi)\chi = -\int \psi (\nabla L^{-1}\chi)
\ee
The adjoint operator of $\nabla L^{-1}$ is $-\nabla L$, given by
\be
  \label{ii55}
  \big(\nabla L^{-1}\big)^+ = -\nabla L
\ee
The Laplacian $\nabla^2$ can be written as $\nabla LL^{-1}\nabla$ and
we obtain in this way:
\ba
  \label{ii56}
  \int \psi^{\ast}\nabla^2\psi &=& 
  \int \psi^{\ast}\nabla LL^{-1}\nabla\psi\\
  &=& -\frac{1}{q}\int 
  \big(\nabla L^{-1}\psi^{\ast}\big)\big(\nabla L^{-1}\psi\big)\nonumber
\ea
An equivalent action to (\ref{ii53}) is:
\be
  \label{ii57}
  \int\limits_{t_1}^{t_2}dt\,\int\limits_{2N}^{2M}
  \big\{i\psi^{\ast}\frac{\partial}{\partial t}\psi - 
  \frac{1}{2mq}\big(\nabla L^{-1}\psi^{\ast}\big)
  \big(\nabla L^{-1}\psi\big) - V\psi^{\ast}\psi\big\}
\ee
We have a Lagrangian that depends on $\psi$, $\dot{\psi}$ and $\nabla
L^{-1}\psi$ as well as on the conjugate expresions.

Let us now assume that we have a Lagrangian which depends on the fields
$\psi$, $\dot{\psi}$ and $\nabla L^{-1}\psi$. The $\psi^{\ast}$ is
considered to be an independent field.  The variation of the action
due to the variation of the fields:
\be
  \label{ii58}
  \psi'=\psi+\delta\psi,\qquad 
  \delta\psi(q^{2M},t_2)=\delta\psi(q^{2N},t_1)=0
\ee
can be written as follows:
\ba
  \label{ii59}
  \delta W &=& \int\limits_{t_1}^{t_2}dt\,\int\limits_{2N}^{2M}
  \Big\{\frac{\partial\c{L}}{\partial \psi}\delta\psi + 
  \frac{\partial\c{L}}{\partial\dot{\psi}}\delta\dot{\psi} + 
  \frac{\partial\c{L}}{\partial\,\nabla L^{-1}\psi}
  \delta\nabla L^{-1}\psi\Big\}\nonumber\\
  &=& \int\limits_{t_1}^{t_2}dt\,\int\limits_{2N}^{2M}
  \Big\{\frac{\partial\c{L}}{\partial \psi}\delta\psi + 
  \frac{\partial\c{L}}{\partial\dot{\psi}}\frac{\partial}{\partial t}
  \delta\psi + \frac{\partial\c{L}}{\partial\,\nabla L^{-1}\psi}
  \nabla L^{-1}\delta\psi\Big\}\\
  &=& \int\limits_{t_1}^{t_2}dt\,\int\limits_{2N}^{2M}
  \Big\{\frac{\partial\c{L}}{\partial \psi} - 
  \frac{\partial}{\partial t}\frac{\partial\c{L}}{\partial\dot{\psi}} - 
 \nabla L\frac{\partial\c{L}}{\partial\,\nabla 
L^{-1}\psi}\Big\}\delta\psi\nonumber
\ea
The last step involves the fact that the variations of the fields
vanish at the boundary. We obtain thus the Euler-Lagrange equation:
\be
  \label{ii60}
  \frac{\partial\c{L}}{\partial \psi} = 
  \frac{\partial}{\partial t}\frac{\partial\c{L}}{\partial\dot{\psi}}
   + \nabla L \frac{\partial\c{L}}{\partial\,\nabla L^{-1}\psi}
\ee
It is easy to verify that eqn. (\ref{ii60}) yields the Schr\"odinger
equation (\ref{ii1}) for the action (\ref{ii57}).

To formulate the Noether theorem we study variations of the field
$\psi$ at the same spacetime point
\be
  \label{ii61}
  \psi'(x)=\psi(x)+\Delta\psi(x)
\ee
which leave the action $W$ invariant
\be
  \label{ii62}
  \delta W = \int dt\,\int \big(\c{L}'-\c{L}\big) = 0
\ee
We expand $\c{L}'$,
\be
  \label{ii63}
  \c{L}' = \c{L} + \frac{\partial\c{L}}{\partial \psi}\Delta\psi + 
   \frac{\partial\c{L}}{\partial\dot{\psi}}\Delta\dot{\psi} + 
   \frac{\partial\c{L}}{\partial\,\nabla L^{-1}\psi}\Delta(\nabla L^{-1}\psi)
\ee
and we find that
\be
  \label{ii64}
  \Delta\c{L} = \c{L}'-\c{L} = \frac{\partial\c{L}}{\partial \psi}\Delta\psi + 
 \frac{\partial\c{L}}{\partial\dot{\psi}}\Delta\dot{\psi} + 
 \frac{\partial\c{L}}{\partial\,\nabla L^{-1}\psi}(\nabla L^{-1}\Delta\psi) =0
\ee
For the first term we insert the Euler-Lagrange equation (\ref{ii60})
and we obtain
\be
  \label{ii65}
  \frac{\partial}{\partial t}\Big\{\frac{\partial\c{L}}
  {\partial\dot{\psi}}\Delta\psi\Big\}+
  \nabla\Big\{\Big(L\frac{\partial\c{L}}{\partial\,
  \nabla L^{-1}\psi}\Big)\big(L^{-1}\Delta\psi\big)\Big\} = 0
\ee
This is Noether's theorem.  If we apply it to the Lagrangian
(\ref{ii57}) and the phase transformation we obtain
\be
  \label{ii66}
  \psi'=e^{i\alpha}\psi,\quad \Delta\psi = i\alpha\psi,\quad 
  \Delta\psi^{\ast}=-i\alpha\psi^{\ast}
\ee
For the ``charge'' density we find from (\ref{ii65})
\be
  \label{ii67}
  -\alpha\rho = \frac{\partial\c{L}}{\partial\dot{\psi}}\Delta\psi+
  \frac{\partial\c{L}}{\partial\dot{\psi}^{\ast}}\Delta\psi^{\ast} = 
  -\alpha \psi^{\ast}\psi
\ee
and for the current
\ba
  \label{ii68}
  &&\hspace{-2cm}\Big(L\frac{\partial\c{L}}{\partial\,
  \nabla L^{-1}\psi}\Big)\big(L^{-1}\Delta\psi\big)+
  \Big(L\frac{\partial\c{L}}{\partial\,\nabla L^{-1}\psi^{\ast}}\Big)
  \big(L^{-1}\Delta\psi^{\ast}\big) =\nonumber\\
  &=& -\frac{i\alpha}{2mq}\big\{(L\nabla L^{-1}\psi^{\ast})
  (L^{-1}\psi)-(L\nabla L^{-1}\psi)(L^{-1}\psi^{\ast})\big\}\nonumber\\
  &=& +\alpha\frac{1}{2im}\big\{(\nabla\psi^{\ast})
  (L^{-1}\psi)-(\nabla\psi)(L^{-1}\psi^{\ast})\big\}\\
  &=& - \alpha\frac{1}{2im}L^{-1}\big\{\psi^{\ast}
  (L\nabla\psi)-(L\nabla\psi^{\ast})\psi\big\}\nonumber\\
  &=& -\alpha j\nonumber
\ea
This agrees with our definition of $j$ in eqn. (\ref{ii4}). There we 
have verified explicitely the continuity equation (\ref{ii3}).

\sect{The $q$-Harmonic Oscillator}

In analogy to~\cite{Lor96} we define a $q$-deformation of the Harmonic
Oscillator with the help of a creation and annihilation operator:
\ba
  \label{iii1}
  a &=& \alpha L^{-2} - i\beta  \nabla L^{-1}   \nn
  a^+ &=& \bar{\alpha} q^{-2} L^2 - i\bar{\beta} \nabla L 
\ea
They satisfy the $q$-commutation relation:
\be
  \label{iii2}
  a a^+ = q^{-2} a^+  a + q^{-2}(1-q^{-2})|\alpha|^2
\ee
To normalise this equation we set:
\be 
 \label{iii3}
 |\alpha| = \frac{q}{\sqrt{1-q^{-2}}}
\ee
Here $\alpha$ and $\beta$ are only up to $q$-factors equal to the
respective constants in~\cite{Lor96}.  \

A general expression for the operators in (\ref{iii1}) can be found:
\ba
  \label{iii4}
   a &=& \alpha L^{-2m} - i\beta L^{-m-1} \nabla L    \nn
   a^+ &=& \bar{\alpha} q^{-2m} L^{2m} - 
   iq^{-m-1}\bar{\beta} \nabla L^{-1} L^{m+1}
\ea
with the commutation relation:
\be
  \label{iii5}
  a a^+ - q^{-2m} a^+ a  = q^{-2m}(1-q^{-2m}) |\alpha|^2
\ee
The Hamiltonian for the $q$-deformed oscillator based on (\ref{iii1}),
(\ref{iii2}) has the following form:
\be
   \label{iii6}
    H = a^+ a = |\alpha|^2 q^{-2} - 
    i \bar{\alpha} \beta \nabla L - 
    i \alpha \bar{\beta} \nabla L^{-1} - q|\beta|^2 \nabla^2
\ee
and the Schr\"odinger equation
\be
  \label{iii7}
  i\partial_t \psi = \left( -q|\beta|^2 \nabla^2 - 
  i\bar{\alpha}\beta \nabla L - i\alpha \bar{\beta} 
  \nabla L^{-1} + q^{-2} |\alpha|^2 \right) \psi
\ee
The Lagrange function to this equation is:
\ba
   \label{iii8} 
    \c{L} &=&  i \psi^* \partial_t \psi - |\beta|^2 
    \left( \nabla L^{-1} \psi^* \right) \left( \nabla L^{-1}\psi \right)  \nn
    && - i\bar{\alpha} \beta 
   \left( \nabla L^{-1} \psi^* \right)\psi +  
    i\alpha \bar{\beta} \psi^* \left( \nabla L^{-1}\psi \right) - 
    |\alpha|^2 q^{-2} \psi^* \psi  
\ea

Now we have a look at the spectrum. We first consider the ground state.  
The Hamiltonian of the $q$-deformed Harmonic Oscillator is a positive operator. 
So we define a ground state by:
\be
  \label{iii9}
   a|0\rangle = 0
\ee
We have a $q$-difference equation:
\be
  \label{iii10}
  \alpha L^{-2} \psi_0(x) = i\beta \nabla L^{-1} \psi_0(x)
\ee
When we consider an Ansatz with a $q$-deformed exponential function:
\be
  e_{q^{-2}}(x) \equiv \sum_{k=0}^\infty \frac{x^k}{(q^{-2};q^{-2})_k}
\ee
we find for $\psi_0(x)$:
\be
  \label{iii11}
  \psi_0(x) = e_{q^{-2}}\left( -i\frac{\lambda \alpha}{q^2 \beta} x \right)
\ee
This function can be seen to be the $q$-Fourier transform of the
gaussian function:
\be
 \label{iii12}  
 f(q^l)  =  q^{-\frac 12(l^2 +l)}c_0
\ee
We consider the $q$-Fourier transformation of (\ref{i63}) and consider the
Ansatz:
\ba
  \label{iii13}
  g(\tau q^{2\nu}) &=& \frac{N_{q}}{\sqrt{2}}\sum_{l=-\infty}^{\infty}
  q^{\nu+l}\Big(f(q^{2l})\cos_q2(\nu+l) \nn\nn
  && +i\tau f(q^{2l+1})\sin_q2(\nu+l)\Big)    
\ea
Inserting the definitions of $\cos_q2(\nu +l)$, $\sin_q2(\nu +l)$ and
noting that the sum over $l$ can be cast into a constant by a
Gauss-summation:
\be
 \label{iii14}
 \sum_{l=-\infty}^{\infty} q^{-2(l-n)^2} c_0 \equiv \tilde{c_0}
\ee
which takes the value (by Jacobi's Triple Product Identity~\cite{Gas90}): 
\be
  \tilde{c_0} = c_0 (q^{-4},-q^{-2},-q^{-2};q^{-4})_\infty
\ee
we find that
\ba
  \label{iii15}
  g(\tau q^{2\nu}) &=& \frac{N_{q}}{\sqrt{2}} \tilde{c_0} q^\nu 
  \sum_{n=0}^\infty (-1)^n \frac{q^{-2n}}{(q^{-2};q^{-2})_{2n}} 
  q^{4\nu n}  \nn
  &&  + i \tau \sum_{n=0}^\infty (-1)^n 
  \frac{q^{-(2n+1)}}{(q^{-2};q^{-2})_{2n+1}} q^{2\nu (2n+1)}\nn
  &=&  \frac{N_{q}}{\sqrt{2}} \tilde{c_0} q^\nu e_{q^{-2}} 
  (i \tau q^{-1}q^{2\nu})  
\ea
The two sums are precisely the real and imaginary part of a $q$-deformed
exponential function. The function $f(q^l)$ is the ground state function
from~\cite{Lor96} (with $\frac{\alpha}{\beta} = 1$) which was calculated 
in the momentum basis.  So we explicitly calculated the $q$-Fourier 
transformation of the ground state in momentum space to the ground state 
in configuration space which is a $q$-deformed exponential function.  
For the odd components we consider the Ansatz:
\ba
  \label{iii16}
  g(\tau q^{2\nu +1})   &=& 
  \frac{N_{q}}{\sqrt{2}}\sum_{l=-\infty}^{\infty}  
  q^{\nu+l}\big( f(q^{2l+1})q\cos_q2(\nu+l+1) \nn\nn
  && +i\tau f(q^{2l})\sin_q2(\nu+l)\big) 
\ea
and find by an analogous calculation:
\be
  \label{iii17}
  g(\tau q^{2\nu +1}) = \frac{N_{q}}{\sqrt{2}}c_0'q^\nu e_{q^{-2}}
  (i\tau q^{-1}q^{2\nu+1})
\ee
The constant $c_0'$ is different from $\tilde{c_0}$:
\be
  c_0' = c_0 (q^{-4},-q^{-4},-1;q^{-4})_\infty
\ee
The functions (\ref{iii15}) and (\ref{iii17}) are the same as
(\ref{iii11}) evaluated in the representation of ~\cite{Hin98}  on 
even and odd lattice points respectively:
\be
  \label{iii18}
  X|\nu, \tau \rangle = -\tau \frac{q^{\nu}q^{- \frac 12}}
  {\lambda}|\nu, \tau \rangle
\ee
Note that in comparison to the Notation in~\cite{Lor96} the constant
 $\frac{\alpha}{\beta}$ here takes the value $q^{\frac{3}{2}}$.

Next we consider the excited states of the $q$-oscillator. We follow
the arguments of~\cite{Lor96} and find that the $q$-deformed Hermite
Polynomials appear in the same way. First we note that:
\be
  \label{iii19}
  a^+|0\rangle = i\frac{1}{q^{\frac 12}\beta} X|0\rangle 
\ee
which is easily verified with the help of  (\ref{iii1}) and (\ref{iii9}). 
To show that the n-particle state can analogously be
expressed by a polynomial in $X$:
\be
 \label{iii20}
 \left( a^+ \right)^n|0\rangle = 
 \left( \frac{1}{\sqrt{2}}\right)^n H_n^{(q)}
 \left( \frac{iX}{\sqrt{2} \beta}\right) |0\rangle
\ee
We calculate the commutation relation between $a^+$ and $\xi$, where
$\xi$ is dimensionless:
\be
  \xi \equiv \frac{i}{\sqrt{2}\beta} X  
\ee
and find:
\be
  \label{iii21}
  a^+ \xi = q^{-2} \xi a^+ - \frac{q^{-\frac{3}{2}}}{\sqrt{2}}
\ee
This leads to the same recursion relation as in~\cite{Lor96} for the
polynomials $H^{(q)}_n(\xi)$, because the recursion relation can be
proven by induction over $n$ with the help of (\ref{iii21}):
\be
  \label{iii22}
  H_{n+1}^{(q)}(\xi) - 
  q^{-\frac 12}q^{-2n}2\xi H_{n}^{(q)}(\xi) + 
  2q^{-n-1}[n]H_{n-1}^{(q)}(\xi) = 0
\ee
We have used here the symmetric $q$-number
$$
  [n] \equiv \frac{q^n-q^{-n}}{q-q^{-1}}
$$

Finally we want to consider another example of the $q$-Fourier transform of 
\cite{Koo92}.
We calculate the $q$-Fourier transform of the step function which will be useful 
for
further applications. We define the step function
\be
     \label{iii23}
          \Theta(q^{2n} - q^{2M}) \equiv   \left\{\begin{array}{rl} 1 & \quad ,n 
\le M\\
                                            0 & \quad ,n > M  
                          \end{array}\right.\nonumber
\ee
Now we calculate the $q$-Fourier transform of this function:
\begin{eqnarray}
  \label{iii24}
            \widetilde{\Theta}( q^{2k} -q^{2M}) &=& N_q \sum_{n=-\infty}^\infty 
q^{2n}
                            \cos_q(q^{2(k+n)} \Theta(q^{2n} - q^{2M})      \nn
                       &=&   N_q \sum_{n=-\infty}^M q^{2n} \cos_q(q^{2(k+n)})     
\end{eqnarray}
The sum can be calculated with the help of the integral (\ref{i46}) and the fact 
that $\cos_q(z)$ can be expressed as the $q$-derivative of $\sin_q(z)$ 
(\ref{i75}). We find using Stokes theorem
\begin{equation}
  \label{iii25}
              \widetilde{\Theta}( q^{2k} -q^{2M}) = N_q q^{-2k} 
\sin_q(q^{2(k+M)})  
\end{equation}
This result can also be obtained by applying the $q$-difference relation 
(\ref{i66})
to the sum in (\ref{iii24}).

To verify our calculation we try to perform the $q$-Fourier transformation in 
the other
direction, too:
\begin{eqnarray}
  \label{iii26}
                  g(q^{2n}) &=& N_q \sum_{k=-\infty}^\infty q^{2k}
                           \cos_q(q^{2(k+n)} \widetilde{\Theta}( q^{2k} -q^{2M}) 
 \nn
                         &=& N_q^2 \sum_{k=-\infty}^\infty \cos_q(q^{2(k+n)}
                                          \sin_q(q^{2(k+M)})
\end{eqnarray}
The sum on the right hand side can be written in terms of $q$-Bessel functions, 
for the notation see \cite{Koo92}:
\begin{equation}
  \label{iii27}
           g(q^{2n}) = q^{M+n}   \sum_{k=-\infty}^\infty q^{-2k} 
J_{-\frac{1}{2}}(q^{2(n-k)})
                                                   J_{\frac{1}{2}}(q^{2(M-k)})
\end{equation}
To calculate this sum we use the summation formula  for $q$-Bessel 
functions from \cite{Koo92} and find
\begin{eqnarray}
  \label{iii28}
       g(q^{2n}) = q^{M+n} \left\{ \begin{array}{c}
                                              
q^{M-3n}\frac{(1,q^{-6};q^{-4})_\infty}{(q^{-2},q^{-4};q^{-4})_\infty}
                                                      
{}_2\Phi_1(q^{-4},q^{-2};q^{-6};q^{-4},q^{-4(n-M)})  \nn
                                                                      \nn
                                               
q^{-(M+n)}\frac{(q^{-4},q^{-2};q^{-4})_\infty}{(q^{-2},q^{-4};q^{-4}_\infty)}
                                                  
{}_2\Phi_1(1,q^{-2};q^{-2};q^{-4},q^{-4(M-n+1)})
                                                      
\end{array}\right.\nonumber
\end{eqnarray}
Because of the finite radius of  convergence of the hypergeometric series which 
is $r=1$ in 
this case \cite{Gas90}, we get different conditions for the variable $n$. Noting 
that  
\begin{equation}
  \label{iii29}
             (1;q^{-4})_\infty = 0
\end{equation}
We find
\begin{equation}
  \label{iii30}
          g(q^{2n}) = \left\{ \begin{array}{rl}
                                         0 & \quad ,n>M  \nn
                                         1 & \quad , n\le M
                                       \end{array}\right. 
\end{equation}
This is exactly the step function we defined in (\ref{iii23}).

\section*{Acknowledgment}
The authors would like to thank Stefan Schraml for useful discussions.
BLC would like to thank the `Fondazione Angelo della Riccia'.

\end{document}